\documentclass[12pt]{article}

\newtheorem{theorem}{Theorem}[section]
\newtheorem{problem}[theorem]{Problem}
\newtheorem{lemma}[theorem]{Lemma}
\newtheorem{proposition}[theorem]{Proposition}
\def\debproof{{\bf Proof: \ }}
\def\finproof{\null\hfill {$\Box$}\bigskip}
\usepackage{latexsym}
\usepackage{graphics}

\usepackage{amsfonts}
\usepackage{times}

\begin{document}
\title{Patterning by genetic networks and modular principle}
\author{S. Genieys \and S. Vakulenko
}                     
%
%
%
%
%
\maketitle
\begin{abstract}
We consider here the morphogenesis (pattern formation)
problem
for some genetic network models.
First, we show that any given spatio-temporal
 pattern
  can be generated
by a genetic network involving a sufficiently large number
of genes.
Moreover,  patterning process can be  performed  by
 an effective algorithm. We also show that
Turing's or Meinhardt's type  reaction-diffusion
models  can be approximated by
  genetic networks.

These results exploit
the fundamental
fact that  the genes form functional units and are organised in blocks
(modular principle).
Due to this modular organisation, the genes  always
are capable to construct any new patterns and even any time sequences
of new patterns from old patterns.
Computer simulations illustrate   analytical results.
\end{abstract}
\section{Introduction}
\label{intro}
In this paper we consider
 pattern formation  problem in the developmental biology.
Mathematical approaches to this problem start with the seminal
work by A. M. Turing \cite{Turing} devoted to pattern formation
 from
a spatially uniform state. Turing's model is a system
of two reaction-diffusion
equations. After \cite{Turing},  similar
phenomenological models were studied
by numerous works
(see \cite{Meinhardt,Murray} for
the review). Computer simulations based on
this mathematical approach  give patterns
 similar to really observed ones
\cite{Meinhardt}. However,
 there is no direct evidence of Turing's patterning in any
 developing organism (\cite{Wolpert}, p.347).
The mathematical models
are often
selected to be mathematically tractable and they do not take
into account actual experimental genetic information.

Moreover, within  the framework of the Turing-Meinhardt
approach some important theoretical questions are
left open. For example, whether there exist  "universal"
mathematical models and
 patterning algorithms
that allow to obtain any, even very  complicated, patterns.
In fact, a difficulty in  using of simple reaction-diffusion
models with polynomial or rational nonlinearities is
that
we have no patterning algorithms.
To obtain a given pattern, first we  choose
a reasonable model (often using
intuitive ideas) and later
we adjust
coefficients or nonlinear terms
by numerical experiments
(an excellent example of this approach is given by
 the book of H. Meinhardt on pigmentation in shells
\cite{Meinhardt2}).

To overcome this algorithmic difficulty
we use genetic circuit models. We are going to show that
they can serve as "universal models", which are capable to generate any
spatio-temporal patterns by algorithms.
The gene circuits were proposed and investigated by many works
\cite{Edw,Glass,Mjol,Rein1,Esp,Rein2,Smolen,Thomas} (for the review see \cite{Smolen})
in order to use available genetic information, to take into
account some fundamental properties of gene interaction
and understand mechanisms of cell gene regulation.

In this paper we investigate the model from \cite{Mjol,Rein1},
which is similar to the well studied Hopfield neural networks.
This model
 describes activation or depression of one gene by another
and have the following form:

\begin{equation}
 \frac{\partial y_i}{\partial t}=R_i\sigma( \sum_{j=1}^m
K_{ij} y_j -\theta_i(x)-\eta_i) -\lambda_i y_i +  d_i \Delta y_i,
\label{genemodel}
\end{equation}
where $m$ is the number of genes included in the circuit,
$y_i(x,t)$ are the concentration of the $i$-th protein,
$\lambda_i$ are the protein decay rates, $R_i$ are some positive
coefficients
and $d_i$ are the protein diffusion coefficients.
We consider (\ref{genemodel}) in some  bounded domain
$\Omega$ with a  boundary $\partial \Omega$.


The real number $K_{ij}$ measures the influence of the $j$-th gene
on the $i$-th one. The assumption
 that gene interactions can be expressed by
a single
real number per pair of genes is a simplification excluding complicated
interactions between three, four and more genes.
Clearly such interactions are possible, however
in this case
the problem becomes
mathematically much more complicated. Since the pair interaction
is capable to produce
any  patterns,
it seems reasonable to restrict our consideration only to such
interaction.

The parameters $\eta_i$ are activation thresholds  and
$\sigma$ is a monotone function
satisfying the following assumptions
\begin{equation}
\sigma \in C^{\infty}({\bf R}), \lim_{z \to -\infty} \sigma(z)=0,
\lim_{z \to +\infty} \sigma(z)=1,
\label{sigma1}
\end{equation}
\begin{equation}
|\frac{d\sigma}{dz}| < C\exp(-c|z|), \quad \sigma^{\prime}(0)=1.
\label{sigma2}
\end{equation}

The well known example is
 $\sigma(z) ={(1 + \tanh z)}/{2}$.


The functions $\theta_i(x)$ are
other activation thresholds
depending on $x$. They can
be interpreted as densities of proteins associated
with the maternal genes.

This model takes into account only three fundamental processes:
(a)  decay of gene products (the term $-\lambda_i y_i$);
(b) exchange of gene products between cells (the term with $\Delta $)
and (c) gene regulation and protein synthesis.
Notice that this model of gene circuit can be considered
as a Hopfield's neural network \cite{Hop} with thresholds
depending on $x$ and where diffusion is taken into account. 
The Hopfield system is the first model
of so-called attractor neural network, both fundamental
and simple.   Analytical
methods for the Hopfield models  were developed
in \cite{Fun,Vak1,Vak5,Vak4}.

Let us fix a function $\sigma$ satisfying (\ref{sigma1}), (\ref{sigma2})
and functions $\theta_i$.
On the contrary,
we  consider
$m, {\bf K},
\lambda_i, d_i,  R_i$ and $\eta_i$ as parameters to be
 adjusted. We  denote
the set of these parameters  by ${\bf P}$ :
\begin{equation}
{\bf P}= \{m, {\bf K}, \eta, {\lambda}, {d}, R \}.
\label{P}
\end{equation}

  Model (\ref{genemodel})
 allows to
use  data on gene regulation (see
\cite{Rein1}, where the
least square approximation of experimental data and
simulated annealing were used to determine the values of
the parameters ${\bf P}$).

In order to study (\ref{genemodel}), many previous works used
numerical simulations. For example,
the work \cite{Rein2} is devoted to
the segmentation in {\sl Drosophila}, in \cite{Esp} the authors
analyse complex patterns occurring under a random choice
of the coefficients $K_{ij}$.

Let us formulate now mathematically our main problem.

\begin{problem}[Universal pattern generation problem] 
Let $T_0>0$ and
$T_0 < T$. Given a function $z(x,t), x \in \Omega, t \in [0,T]$ and
a positive number
$\epsilon$, to find the parameters ${\bf P}$
such that the solution  of
system (\ref{genemodel})  with initial conditions  $y_j =0$ satisfies
\begin{equation}
 \sup_{x,t} |z(x,t) - y_1(x,t)| < \epsilon, \quad x \in \Omega,
\quad t \in [T_0, T].
\end{equation}
\end{problem}

Let us consider a biological interpretation
of this mathematical formulation.
We assume that
the cell states depend on expression of
some genes; we can thus identify observed patterns
of cell differentiation with gene expression patterns.
Let us  consider  these expression patterns as continuous functions
 $z(x, t)$, where $x \in \Omega, t  \in [0, T]$
 and $T > 0$.  For example, we can assume that
$z \in [0,1]$ and if $z $ is close to $1$, the gene is expressed
otherwise the gene is not expressed.

We consider  gene circuits  including a single "output"
(structural)
gene   $y_1=Y_{out}$ and $m_1$ "hidden"  (regulating) genes.
The output gene  can change the cell states  and therefore
 can predetermine an output pattern
$z$. The hidden genes do not influence directly the cell states,
 they are involved only in an internal cellular gene regulation.

Notice that given pattern $z$ can depend on time $t$.
This fact is important since real biological structures
are usually dynamical.
For stationary patterns $z$ (independent of time)
the solution of pattern problem is simple and follows from
the well known results on neural networks
(see Sect. ~\ref{sec-aux}).

The main results  can be described as follows.

{\bf A)}
 We show that, roughly speaking,  any pattern formation process based
on a reaction-diffusion model can be performed
as well by a
 genetic network, with a sufficiently large number of the genes.
For each reaction-diffusion model
one can find an approximating gene network, with the almost same
pattern formation capacity.
This result  justifies, to some extent,   Turing-Meinhardt's
 models from a genetic point of view.
Indeed, these models
can be  considered as  gene
circuits.

{\bf B)}
The second result asserts that, under
natural conditions on maternal genes
densities $\theta_i$,
the universal pattern generation
problem always has a solution.
Moreover, there is a constructive and numerically effective
algorithm that allows us to find a circuit
generating a given pattern.

Notice that this result
is also valid in the absence of diffusion.
Indeed, in our approach, spatial signalling is not provided
by the diffusion process, but by space-depending thresholds $\theta_i$.

Our conditions on maternal gene concentration
are necessary and sufficient:
if they do not hold,
 it is not possible to approximate any patterns within an
arbitrarily small error. On the contrary,
if they are valid, it  is possible.
If we deal with one-dimensional case (for example, we consider
a differentiation along anterior-posterior axis), then
our conditions mean existence of a
morphogene  gradient along this axis. For {\sl
Drosophila} this morphogene is {\sl bicoid}.
So, we show that
a simple bicoid gradient is capable to produce any chain of
complicated time transformations  leading to complex spatial
 one-dimensional patterns.
 This result is in an agreement with
 biological observations \cite{Alberts,Wolpert}.
To create any two-dimensional patterns, we need at least two independent
gradients, along anterior-posterior
and dorso-ventral axes.

{\bf C)}  We show  that the modular
organisation and sigmoidal interaction
are  effective tools to  form complex hierarchical
patterns. 

Indeed, we show that new, more refined structures,
can be obtained by using of  previous old structures.
Also we illustrate that existence of an old structure make it
easier to produce a new complex one. 
This property might help to understand the 
usual idea
"morphogenesis repeats evolution"   \cite{Rid}, see
 Sect. ~\ref{sec-prog} and ~\ref{sec-ccl}.

The paper is organised as follows. In the next section
we explain main biological and mathematical ideas beyond
these results, in particular, we
find connections with multilayered
network theory and  the Hopfield model.
In Sect. ~\ref{sec-RD} we describe the connection
between the reaction-diffusion models and gene circuit
systems.  In Sect. ~\ref{sec-prog}  we formulate mathematically
pattern formation problem and describe main ideas of
patterning algorithms.
Sect. ~\ref{sec-num} presents computer simulations illustrating our analytical
results. In particular,
as an illustration, we approximate numerically,
by gene circuit, a reaction-diffusion
system for  pigmentation of sea shells proposed by \cite{Meinhardt2}.
Sect. ~\ref{sec-ccl} contains a
discussion and concluding remarks.

All  complicated and tedious mathematical  details
can be found in the Appendix.

\section{Main mathematical instruments}
\label{sec-aux}

In this section we remind main ideas and
results of neural network theory
important below. To simplify our statement, we
omit some non-essential mathematical details
(for details, see
\cite{Vak5}).

\subsection{Multilayered neural networks}

The neural networks usually consist of a large number of neurons.
Each neuron
is connected to other neurons by directed links with their
associated  weights. After absorbing the inputs, each
neuron produces its activation as an output signal to other
neurons. Each neuron sends a single signal to several neurons at
the  time.   Typical problems which may be solved by such nets
are pattern classification, storing patterns and optimal control
problems (see \cite{Fried}).

The simplest example is a single-layer network
having one layer of weights.   The network consists of
$n$ input neurons $X_j$, $j=1,\dots,n$ and an output neuron $Y$.
Each $X_j$ is connected to  $Y$ with an associated weight
$w_j$. The output $Y$ is given by

\begin{equation}
  Y=\sigma(\sum_{j=1}^n w_j X_j-h),
\label{s-net}
\end{equation}
where $h$ is a threshold and $\sigma$ is a strictly monotone function
satisfying  (\ref{sigma1}) and (\ref{sigma2}).

Network (\ref{s-net}) can solve only simple classification problems.
More powerful, a multilayer neural network with $p$ layers consists
of one layer of input neurons,  an output neuron and $(p-1)$
hidden layers. For $p=2$, the corresponding equations
can be written for instance as
\begin{equation}
  Y=\sigma(\sum_{k=1}^m B_k z_k-h),
\label{s1-net}
\end{equation}
\begin{equation}
  z_k=\sigma(\sum_{j=1}^n A_{kj} q_j-\eta_k),
\label{s2-net}
\end{equation}

where $q_j$ are states of the input neurons.
The remarkable property of this network playing the key role in this paper
is that
any input-output map
of the form  $(q_1, q_2,\dots,q_n) \to F(q_1, q_2,\dots, q_n)$, where
$F$ is a continuous function,
can be approximated by a network (\ref{s1-net})--(\ref{s2-net})
with a sufficiently large $m$ and
appropriate weights $A_{kj}$ and $B_k$.

We shall use below, for brevity, notation
\begin{equation}
{\bf A}_j q =\sum_{l=1}^m A_{jl} q_l.
\label{sum}
\end{equation}
Since $\sigma$ is monotone, the assertion that
 network (\ref{s1-net})--(\ref{s2-net}) approximates any output, 
can be reformulated
as follows: by the quantity
$\Psi =\sum_{k=1}^m B_k z_k$ we can approximate any function,
within  an arbitrarily small error.
This fact
results from the following well known
assertion (see, for example,  the works \cite{Barron,Fun,Hor,Vak5}).

Let us consider function $\Psi(q, {\bf A},
{ B}, \eta)$ of vector argument $q=(q_1,\dots, q_m)\in{\bf R}^m$
 depending
on the following parameters:
the number $m > 1$, an $m \times n$ matrix ${\bf A}$,
and the vectors $B$ and
 $\eta \in {\bf R}^m $. This function is defined
by
\begin{equation}
\Psi(q, {\bf A}, { B}, m, \eta)= \sum_{k=1}^m B_k z_k =
\sum_{k=1}^m B_{k} \sigma({\bf A}_k  q - \eta_k).
\label{Psi}
\end{equation}

We consider this function in a bounded ball $\Omega_R$,
consisting of vectors $q$ such that
$|q|^2 =q_1^2 + \dots + q_m^2 < R^2$.

\begin{lemma}[Approximation Lemma]
If $\sigma$ is monotone
and satisfies con\-di\-tions (\ref{sigma1})--(\ref{sigma2}), then
for any continuous  function $Q(q)$ defined
in the ball $\Omega_R$
and for
any positive number $\varepsilon$,
there exist a  number $m \ge n$, matrices
${\bf A}$ and $m$-vector $\eta$ and $B$ such that
\begin{equation}
   |Q(q) -\Psi(q, {\bf A}, { B}, m, \eta)|
< \varepsilon, \quad q \in \Omega_R.
\label{approx}
\end{equation}

In other words, given pattern $Q$ and  $\epsilon > 0$,
we can always find  weights $\bf A$ and $ B$
such that the output of network  (\ref{s1-net})--(\ref{s2-net})
approximates this pattern, up to precision $\epsilon$.

\label{approxlemma}
\end{lemma}

Approximation
(\ref{approx})
can be obtained by the multilayered network
theory (see \cite{Fun,Hech,Hor}) or by an application of
wavelet type extensions \cite{Vak5}. For wavelet theory see \cite{Meyer}.
It is well known that the approximations
by (\ref{Psi}) are
numerically effective as the dimension $n$ of the vector $q$
increases. We know constructive  algorithms allowing to adjust
 the parameters ${\bf A}, { B}, m, \theta$
(see \cite{Barron} and references therein).

Notice that the number $m$ of the coefficients $B$
depends polynomially
on  "complexity" of the pattern $Q$ and the precision $\epsilon$.
More "complex" the pattern, greater $m$. This complexity
can be measured by the integral  \cite{Barron}
\begin{equation}
   \mathrm{Comp}(Q) =\int |\omega| |\hat Q(\omega)| d \omega,
\label{comp}
\end{equation}

where $\hat Q(\omega)$ is the Fourier transform
of the function $Q$. This means that more oscillating
functions (patterns) have larger complexities.

\subsection{Large time behaviour of the Hopfield networks}

The Hopfield network \cite{Hop} is a system of coupled oscillators
defined by  the differential equations
\begin{equation}
 \frac{d y_i}{d t}=R_i\sigma( \sum_{j=1}^m
K_{ij} y_j -\eta_i) - y_i.
\label{Hopfmodel}
\end{equation}

Here $y_i$ are neuron states depending on time, $K_{ij}$ is a matrix
determining a neuron interaction (synaptic matrix), $\eta_i$
are thresholds. Genetic model (\ref{genemodel})
can be considered
as a generalisation of the Hopfield system such that
the neuron states and the thresholds depend on a space variable $x$
and the diffusion is taken into account.

We are going to apply some   methods developed
 to investigate attractors of the Hopfield neural networks.
We recall that
dynamics (\ref{Hopfmodel}) is
dissipative and thus an  attractor always exists
(on the attractor theory see publications \cite{Cu,Ha,Il,La,Ru,T}
among many others).

Dynamics of (\ref{Hopfmodel}) sharply depends
on the synaptic matrix $\bf K$. If
the matrix $\bf K$ is symmetric, the attractor usually consists
of many equilibria. Such stable large time behaviour
can be applied to the pattern recognition and associative memory
problems.

The large time behaviour  of $y$ can become very complex
if $\bf K$ is non-symme\-tric. For instance, depending on  
$\bf K$,
neuron states can form complicated coherent structures
that evolve periodically or even chaotically in time.
These coherent patterns can be described as follows.

Any  (symmetrical or non-symmetrical) $m \times m$ matrix $\bf K$
of rank $n$ can be represented as a product
of two matrices $\bf A$ and $\bf B$, i.e.,
\begin{equation}
  {\bf  K}={\bf A }
{\bf B},
\label{Hopfield}
\end{equation}
where $\bf A$ has size $m \times n$ and
$\bf B$ has size $n \times m$,

Let us introduce the new variables
\begin{equation}
q_l(t)=\sum_{j=1}^m B_{lj} y_j(t)={\bf B}_l y(t),
\label{q}
\end{equation}
where $l=1,2,\dots,n$.

The dynamical equations for $q$ have the
following form
\begin{equation}
\frac{dq_l}{dt}= -q_l +
\Psi_l(q, {\bf A},{\bf B},m, \eta),
 \label{eq-q}
\end{equation}
where $\Psi_l$ are defined by equations similar to (\ref{Psi}).
Time evolution of the new variables  $q_l$
controls the dynamics
of all the neuron states $y_i$.
Indeed,  we have
\begin{equation}
\frac{dy_i}{dt}= - y_i +
\sigma({\bf A}_i q  -\eta_i).
\label{y-q}
\end{equation}

The functions $y(t)$ can be expressed through $q(t)$ in a simple way
by linear equations (\ref{y-q}).




Below we will use
 new control parameters $\bf P$, we denote
${\bf P}=\{n, m, {\bf A}, {\bf B}, \eta \}$ fixing $R_i=1, \lambda_i=1,
d_i=0$.

Let us formulate now the following assertion (analogous
to the results of  \cite{Vak5,Vak4}) describing
the complexity of time behaviour of the  circuits.


\begin{lemma}
\label{control-hopfield}
By the network parameters $\bf P$,  dynamics (\ref{eq-q})
can be specified with\-in an arbitrarily small error.
More precisely, for any $n$, any given continuous functions
$Q_l(q)$ defined
on bounded domain $\Omega$,
and for any $\epsilon > 0$,
we can  choose  parameters
${\bf P}$
such that
\begin{equation}
\label{approx-syst-dyn}
|Q_l + \lambda q_l  -\Psi_l(q,{\bf A}, {\bf B},m, \eta)|
< \epsilon, \quad q \in \Omega, \quad l=1,2,\dots,n.
\end{equation}

Therefore,  any structurally stable dynamics can be generated
by system
(\ref{Hopfmodel}).
\end{lemma}

This
result shows that the  variables
$q_j$ can exhibit complicated dynamics, periodical or chaotical.
In particular, any kind of  stable chaos
can occur in the dynamics of our systems, for example,
the Smale horseshoes, Anosov
flows, the Ruelle-Takens-Newhouse chaos, etc.
\cite{Dyn,NewR,Nit,RuT,Sm,Vi}.

In general, greater the neuron number $m$, more complex this time dynamics.
Thus, the neuron states $y_j$ also can demonstrate a complicated dynamics
however, if $n << m$,
 all the $m$ neuron states are strongly
correlated since they can be defined through
a relatively small number of the hidden variables.

For a proof of  Lemma ~\ref{control-hopfield} see
\cite{Vak5}.

In the next section we shall show
 that the gene networks
can  simulate, in a sense,  any reaction-diffusion
systems.

Notice that some fundamental and simple biological principles are
beyond the mathematics. The genes are organised
in blocks. The local cell differentiation and growth
processes are governed by a collective action of these
blocks.

\section{Approximation of reaction-diffusion
sys\-tems by\\ ge\-ne net\-works}

\label{sec-RD}
We consider, for simplicity,  the case of two component reaction-diffusion
systems

\begin{equation}
\frac{\partial u}{\partial t}=d_1 \Delta u + f(u,v),
\label{RD-1}
\end{equation}
\begin{equation}
\frac{\partial v}{\partial t}=d_2 \Delta v + g(u,v).
\label{RD-2}
\end{equation}

The phenomenological
approach based on (\ref{RD-1})-(\ref{RD-2})  gives
excellent results for
some  pattern formation problems (for example
such as shell pigmentation \cite{Meinhardt2}), where
  nonlinearities can have the following typical form
 \cite{Meinhardt2}
\begin{equation}
f=f_M(u,v)=\alpha v (\frac{u^2}{1+\alpha_1 u^2}+ \beta_1) -\kappa_1 u,
\label{fMeinh}
\end{equation}

\begin{equation}
g=g_M(u,v)=\beta_2 -
\alpha v (\frac{u^2}{1+ \alpha_1 u^2}+ \beta_1) - \kappa_2 v.
\label{gMeinh}
\end{equation}

We suppose here  that all constants $\alpha, \beta_i, \alpha_1,
\kappa_k$ are
positive.

In these equations, $u$ and $v$ are unknown functions of
the space variables
$x=(x_1, x_2, x_3)$ defined in a bounded domain  $\Omega$.

 System (\ref{RD-1})--(\ref{RD-2}) must be complemented
by standard initial and boundary conditions.

Suppose the system of
equations that governs patterning is
 two-component system (\ref{RD-1})--(\ref{RD-2}),
where nonlinearities $f$ and $g$ are
continuous functions.
The general multi-component case can be
studied in a similar way.
Assume
solutions of (\ref{RD-1})--(\ref{RD-2}) remain globally bounded,
i.e., for some positive constants $C_i$ we have the estimate
\begin{equation}
|u(x,t)| < C_1, \quad |v(x,t) | < C_2,
\label{bound}
\end{equation}
for all $t > 0$, if it holds for $t=0$.
Let us define the domain $D_{C_1, C_2}$ as follows:
\begin{equation}
D_{C_1, C_2}=\{ (u,v): \quad  0 \le u < C_1, \ 0 \le v < C_2 \}.
\label{rectangle}
\end{equation}
We suppose  that initial condition
 belongs to $D_{C_1, C_2}$  for each $x$.

Our goal is to show that, for
a given reaction-diffusion system (\ref{RD-1})--(\ref{RD-2})
 we can always find
  an "$\epsilon$- equivalent" circuit (\ref{genemodel}).
Namely, for this equivalent circuit
there exists a smooth map $b(y) : (y_1, y_2,\dots, y_m)
\to (u, v)$ transforming the gene concentrations
to the reagent concentrations and such that  time evolution
of $u, v$ is defined by a new reaction -diffusion system with nonlinearities
$\Phi_1(u,v), \Phi_2(u,v)$,  $\epsilon$- close
to nonlinearities $f(u,v), g(u,v)$.   
Roughly speaking we can say that any reaction -diffusion system
can be realized as a gene circuit.

To this end, we use {\it Modular Principle}.
Let us consider a system (\ref{genemodel}) having  a special block
structure.
Namely,
we assume that there exist two kinds of the genes. We  denote these
groups of the genes
by $y$ and $z$,
where vector $y(x,t)$ contains $m_1$ components
and  $z(x,t)$ contains $m_2$ components. Naturally,
$m=m_1 +m_2$.
We consider  system (\ref{genemodel}) of the special form

\begin{equation}
\frac{\partial y_i}{\partial t}=
\sigma({\bf K}^{yy}_i y + {\bf K}^{yz}_i z - \theta_i)
  +  d_1 \Delta y_i,
\label{gen-RD2-1}
\end{equation}

\begin{equation}
\frac{\partial z_i}{\partial t}=
\sigma({\bf K}^{zy}_i y + {\bf K}^{zz}_i z -\bar \theta_i)
  +  d_2 \Delta z_i.
\label{gen-RD2-2}
\end{equation}

Here we use notation (\ref{sum})
and matrices ${\bf K}^{yy}$, ${\bf K}^{zz}$, ${\bf K}^{zy}$
and ${\bf K}^{yz}$ describe interactions between different groups of
the genes.

In general, these interactions are not symmetric, i.e.,
${\bf K}^{yz}$ is not equal to the transpose of ${\bf K}^{zy} $.

The
coefficients $d_1$ and $d_2$ coincide with the
diffusion coefficients in
equations (\ref{RD-1}) and (\ref{RD-2}).

We choose  the entries of
the matrices ${\bf K}^{yy}$, ${\bf K}^{zz}$, ${\bf K}^{zy}$
and ${\bf K}^{yz}$ as follows:
\begin{equation}
K_{ij}^{yy} = a_i b_j, \quad K_{ij}^{yz} = \gamma_i \bar b_j,
\label{defK1}
\end{equation}
and
\begin{equation}
K_{ij}^{zy} = \bar \gamma_i b_j, \quad K_{ij}^{zz} = \bar a_i \bar b_j
\label{defK2},
\end{equation}
where $a_i, \bar a_i, \gamma_i, \bar \gamma_i, b_i, \bar b_i$
are unknown coefficients.

Let us define "collective variables"
\begin{equation}
u=\sum_{i=1}^{m_1} b_i y_i, \quad v=\sum_{i=1}^{m_2} \bar b_i z_i.
\label{defuv}
\end{equation}

After some calculations ( see the Appendix, part 1)
we obtain
\begin{equation}
\frac{\partial u}{\partial t}=d_1 \Delta u + \Phi_1(u,v),
\label{u}
\end{equation}
and
\begin{equation}
\frac{\partial v}{\partial t}=d_2 \Delta v + \Phi_2(u,v),
\label{v}
\end{equation}
where
\begin{equation}
\Phi_1(u,v)=
\sum_{i=1}^{m_1} b_i \sigma(a_i u + \gamma_i v
- \theta_i),
\label{Phi1}
\end{equation}
\begin{equation}
\Phi_2(u,v)
=\sum_{i=1}^{m_2} \bar b_i \sigma(\bar a_i v + \bar \gamma_i u
- \bar \theta_i).
\label{Phi2}
\end{equation}

Applying  Lemma ~\ref{approxlemma}, we notice
that for any $\epsilon > 0$ there exist
 numbers $m_1$, $m_2$, vectors
$a,b, \bar a, \bar b, \gamma, \bar \gamma$ and $\theta, \bar \theta$
such that
\begin{equation}
|\Phi_1(u,v) - f(u,v)| < \epsilon,
\quad
|\Phi_2(u,v)
  - g(u,v)|  <
\epsilon
\label{approxg}
\end{equation}
for all $u,v$ from some bounded domain.

This proves the  main result of this section:

\begin{proposition}

Con\-si\-der problem (\ref{RD-1})--(\ref{RD-2}) 
whose solutions remain in a domain $D_{C_1, C_2}$.
Then, if functions $f$,
 $g$ are continuous, for any $\epsilon > 0$,
 there exist such a system (\ref{genemodel})
 with a sufficiently large number $m$ and
  coefficients $r=(r_1, r_2, \dots, r_m)$ and
 $s=(s_1, s_2, \dots, s_m)$  such that
 the functions
\begin{equation}
    u= r y=\sum_{i=1}^m r_i y_i,
\quad
 v= s y=\sum_{i=1}^m s_i y_i
\label{combi-uv}
\end{equation}
satisfy the system
\begin{equation}
 u_t= d_1 \Delta  u + \tilde f( u,  v)
\end{equation}
\begin{equation}
v_t= d_2 \Delta v  + \tilde g( u,  v),
\end{equation}
where
\begin{equation}
  |f(u,v) -\tilde f(u,v) | < \epsilon,
\end{equation}
\begin{equation}
|g(u,v) -\tilde g(u,v) |< \epsilon
\end{equation}
for $(u, v) \in D_{C_1, C_2}$.
\end{proposition}

Therefore,  any reaction-diffusion patterning processes on
a bounded time interval $[0,T]$
can  be performed as well by
genetic networks. In other words, the pattern capacity
of the  gene circuits on bounded time intervals are not less than
the pattern capacity
of reaction-diffusion systems.

To conclude this section, let us notice
 that an inverse problem, namely an approximation
of a neural network by a reaction-diffusion system
has been considered in \cite{Edw2} and \cite{Vak5}.

\section{Programming of
spatio-temporal patterns by
 gene circuit models}
\label{sec-prog}

In this section
we  state an analytical
algorithm resolving the following problem:  given spatio-temporal pattern,
to find a gene circuit generating this pattern.
We show that this problem can be
solved even without diffusion ($d_i=0$).
In our approach the space signalling is provided by space-depending
activation thresholds.
It is important
from the biological point of view since the molecular transport
is often performed by non-diffusional mechanisms \cite{Alberts}.
For time discrete networks, similar results
were obtained in \cite{Gr77}.

Beside multilayered network theory (Lemma ~\ref{approxlemma}) we also use
the following result.

\begin{theorem}
[Superposition Theorem]
\label{superp}
Let us consider a family $\cal F$ of gene circuits (\ref{genemodel})
with the parameters
${\bf P}^1, \ {\bf P}^2, \dots, {\bf P}^p$, where
the functions $\theta_i$ are fixed and identical for
all the circuits. Assume
 these networks generate
the output patterns
$Y^1=y_1^1(x,t), Y^2=y_1^2(x,t), \dots, Y^p=y_1^p(x,t)$.

Then, for any $\epsilon > 0$ and for any continuous positive function
$F(u_1,\dots, u_p)$, there is a network
(\ref{genemodel}) generating an output pattern
$y_1=Y$ such that
\begin{equation}
\label{out}
 |F(Y^1(x,t), Y^2(x,t), \dots, Y^p(x,t))-
  Y(x,t)| < \epsilon.
\end{equation}
\end{theorem}

This result can be interpreted as a
 {\it Superposition Principle}. If  given circuits are capable
to produce patterns $Y^1, Y^2,\dots, Y^p$, for any function $F(u_1,\dots,u_p)$
there is a new circuit, which
can approximate the pattern $z$ of the form $z=F(Y^1,\dots,Y^p)$, in other
words, "superposition by $F$"
of these previous patterns.
This result also has interesting biological corollaries; we discuss
it in  Sect. ~\ref{sec-ccl}.

Let us describe first the outline of  the proof.
The proof is based on {\it Modular Principle}.
We suppose that  an unknown interaction matrix ${\bf K}$
of the network can be decomposed  in  blocks. Some blocks contain
 the known matrices ${\bf K}^s$ corresponding to  $s$-th network
of  given network family. An additional block determines
an interaction between  new genes and the  genes involved in
the networks of
 the family $\cal F$.
This structure
allows us to apply  the approximation results of
the multilayered network theory
\cite{Barron,Cyb,Hech,Hor} (see Lemma ~\ref{approxlemma}).
This assumption about the structure of the matrix
$\bf K$ also is in agreement
with contemporary ideas in molecular biology \cite{Hop2,Coll}.
The proof (which,
by {\it Modular Principle}, is quite straightforward)
can be found in the Appendix.

Since the basic element of the proof of   Superposition Principle
is Lemma ~\ref{approxlemma}, and the proof
of this Lemma gives us an algorithm,
 therefore we obtain a complicated but quite
constructive algorithm resolving the patterning problem. Moreover,
we can estimate the number of the genes $N(z)$  involved
in patterning process as a function of the pattern complexity
defined by (\ref{comp}).  Namely, using  the results of
the work \cite{Barron}, we find that $
 N(z)$ depends polynomially on  $\mathrm{Comp}(z|u)$,
where $\mathrm{Comp}(z|u)$ is a conditional complexity of $z$ respectively
given patterns $u$.  To explain this relation and
its biological  meaning, let us consider
a simple example.

Suppose our problem is to construct a periodic one-dimensional
pattern $z(x)=\sin k x$,
where $k$ is  a large number. Our target pattern
therefore is sharply oscillating.
Moreover, we have no stored
(old) patterns $u_i$ and thus $\mathrm{Comp}(z|u)=\mathrm{Comp}(z)$ is proportional to
$k$.
In this case,  to resolve the pattern approximation
problem, the network have to involve many
genes.

Assume now that there are  old patterns $u_i$ and, in particular,
the patterns of the form $\sin k_0 x, \ \cos k_0 x$, where
$k_0 < k$ but $k_0 >> 1$.  In this case the function
$z$ can be expressed through
$u_i$ as a polynom of degree $P=k/k_0$. Thus
$\mathrm{Comp}(z|u)$ is
  much less  $\mathrm{Comp}(z) $ for large
$k_0$ and $k$.

Roughly speaking,
a complex target pattern may be  simple
respectively to another complex pattern.
We discuss a biological interpretation of this property  in Sect. 
~\ref{sec-ccl}.

 Using  Theorem ~\ref{superp}, we can resolve now the pattern programming
problem.
 Suppose  the functions $\theta_i(x)$ possess the following property.
 They can be considered  as  "coordinates" in the domain
 $\Omega$, i.e.,
 there  exist continuous functions $g_i$  such that
\begin{equation}
\label{gy}
x_i=g_i(\theta_1(x), \theta_2(x), \dots, \theta_m(x)),
 \quad x \in \Omega, \ i=1,\dots,n.
\end{equation}
  This condition holds, for example, if $m=n$ and
for each $i$, the function $\theta_i(x)$ is a strictly monotone function
of only one variable $x_i$. A biological example can be given by
the distribution of maternal genes
in Drosophila \cite{Wolpert}.

Let us prove first an auxiliary mathematical result.

\begin{lemma}
\label{changevar}
Suppose that condition (\ref{gy}) holds.
 Then  any con\-ti\-nu\-ous function\\
$F(x_1,\dots,x_n, t)$ can be represented as a function
of $n+1$ variables
\begin{equation}
\bar Y_1=\sigma(\theta_1(x))(1-\exp(-\gamma t)),
\end{equation}
\begin{equation}
 Y_i=\sigma(\theta_i(x))(1-\exp(-\kappa t)),
\end{equation}
for $i=1,\dots,n$, where $\kappa$ and $\gamma$ are two different
positive constants.
\end{lemma}

\debproof
To prove this lemma, let us observe that
\begin{equation}
\log Y_1 - \log \bar Y_1=f(t),
\end{equation}
where $f(t)$ is a strictly monotone function of $t$.
Therefore, $t$ can be written as a function of $Y_1$
and $\bar Y_1$.  Then any $\theta_i(x)$ can be
presented as a function of $Y_i,Y_1$ and $\bar Y_1$.
Using (\ref{gy}), one obtains that each $x_i$
is a function of the variables $Y_s, \ s=1,2,\dots,n$ and $\bar Y_1$.
The lemma is proved.
\finproof

Let us formulate the main result of this work. This result
means that any patterning process can be realized by a gene
circuit.

\begin{theorem}
\label{main}
Suppose that condition (\ref{gy}) holds. Then for any continuous
positive $z(x,t), x \in \Omega, \ t \in [0, T]$, any
positive $T_0 < T$ and  $\epsilon$ there is a
 system (\ref{genemodel})
such that the solution of this system satisfies
the estimate
\begin{equation}
 \label{gy1}
 |z(x,t) - y_1(x,t)| < \epsilon
  \quad x \in \Omega, \quad t \in [T_0, T].
  \end{equation}
\end{theorem}

Before start to prove Theorem, let us notice that
in the case $d_i=0$ (diffusion coefficients vanish)
 condition  (\ref{gy})  is actually
necessary in order to resolve any patterning problem. In other
words, if it does not hold, there is
a pattern, which cannot be $\epsilon$-approximated
for any $\epsilon$. Indeed, if $d_i=0$,
solutions of (\ref{genemodel}) are vector function $y$ of
variables $t$ and $\theta_i$. If  (\ref{gy}) does not hold, for some $s$
the pattern $y_1(x,t)=x_s$ cannot be $\epsilon$-approximated
for any $\epsilon$.
For $d_i \ne 0$  condition  (\ref{gy})  can be replaced
by a weaker one but we will not consider this question here.

\debproof

Theorem ~\ref{main} results  from Lemma ~\ref{changevar} and 
Theorem ~\ref{superp}.
We take a network generating $\bar Y_1$, $ Y_1$, $ Y_2$,
$\dots$, $ Y_n$. This network has the following structure:
\begin{equation}
\frac{ \partial \bar y_1}{\partial t} = -\gamma \bar y_1
+\gamma \sigma (\theta_1),
\end{equation}
\begin{equation}
\frac{ \partial  y_i}{\partial t} = -\kappa  y_i
+\kappa \sigma (\theta_i),
\end{equation}
for $i=1,\dots,n$. We observe now that $y_i=Y_i(x,t)$ and $\bar y_1=\bar Y_1(x,t)$.
This completes the proof.
\finproof

\section{Computer simulations}

\label{sec-num}

We first illustrate the results of Sect. ~\ref{sec-RD}: we approximate the
reaction-diffusion system (\ref{RD-1})--(\ref{gMeinh}) by a gene network.
For this we approximate functions (\ref{fMeinh})
and (\ref{gMeinh})
by sigmoidal functions in order to satisfy inequalities (\ref{approxg}).

This is a problem of nonlinear approximation since functions
$\Phi_1$ and $\Phi_2$ depend linearly on $b_i$ and $\bar b_i$ but
nonlinearly on $a_i$, $\gamma_i$, $\theta_i$ and
$\bar a_i$, $\bar{\gamma}_i$, $\bar{\theta}_i$. Coefficients
$b_i$ and $\bar b_i$ can be  calculated by the classical
least square method, but other coefficients have to be determined
in a proper way. For instance, we could choose these coefficients
randomly and select
the best values (or a satisfying value), however this is
usually too long when the search space is large.
Here this random method cannot be used due to the fact
that the functions (\ref{fMeinh}) and (\ref{gMeinh}) depend on two variables.
If we approximate a function
of a single variable, this simple
random method can be useful and we apply it below.

To make this nonlinear approximation, we use an
iterative approach proposed by Jones \cite{Jones} (see also Barron
\cite{Barron}). Jones'
result is an iterative version of Approximation Lemma ~\ref{approxlemma}: under
the  conditions of this lemma  one can find  a sequence
of approximating functions
$\left(\Psi(q,A,B,\eta,m)\right)_{m\in\mathbb{N}}$ (denoted
$\left(\Psi_m\right)_{m\in\mathbb{N}}$ for brevity) satisfying
\begin{equation}
|Q(q)-\Psi_m|^2\leq C/m,
\end{equation}
where $C$ is a constant depending on $Q$. As already mentioned in 
Sect. ~\ref{sec-aux},
Barron \cite{Barron} has related $C$
to the Fourier transform of $Q$. For more oscillating
 functions $Q$,    the constant $C$ is greater.
This constant can be considered as a measure of  complexity of
$Q$.

\begin{figure}[!hb]
\resizebox{\textwidth}{!}{
\includegraphics{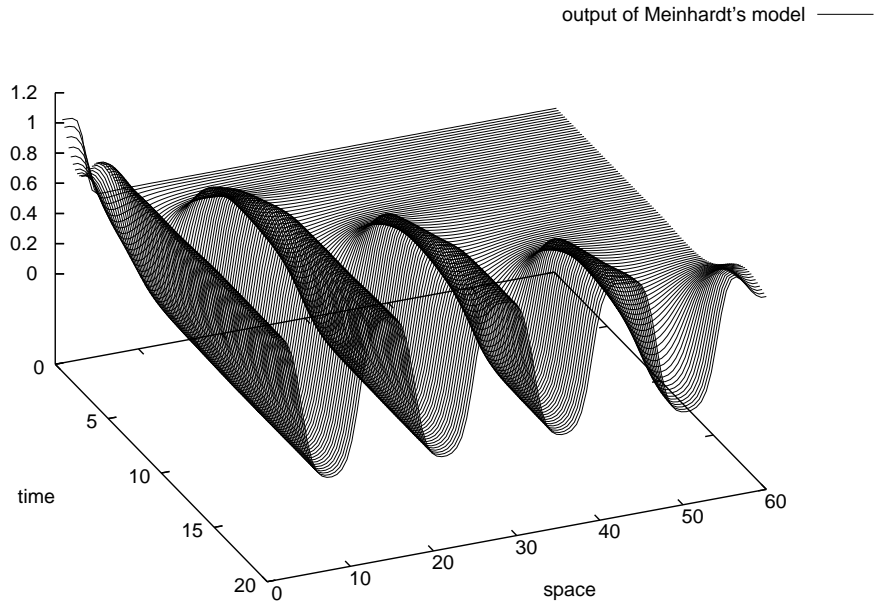}
}
\caption{Solution of Meinhardt's model.}
\label{Meintarget}
\end{figure}

\begin{figure}[!hb]
\resizebox{\textwidth}{!}{
\includegraphics{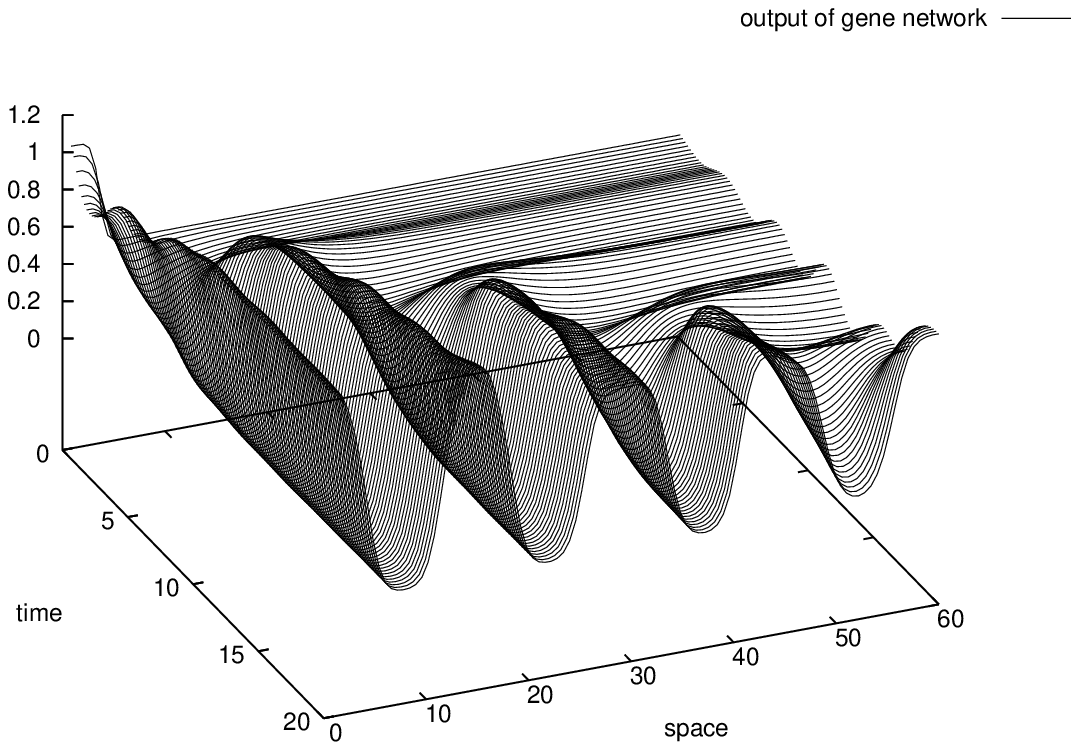}
}
\caption{Approximation of Meinhardt's model by a gene network.}
\label{Meinapp}
\end{figure}

Jones' sequence is defined as follows:
$\Psi_1=B_1\sigma(A_1q-\eta_1)$, where $B_1$, $A_1$, $\eta_1$ give an
almost minimal value of $|Q(q)-\Psi_1|$.
Then, for any $m\geq 2$, $\Psi_m=\alpha_m\Psi_{m-1}+B_m
\sigma(A_mq-\eta_m)$, where $\alpha_m$, $B_m$, $A_m$ and $\eta_m$
give an almost minimal value of
$|Q(q)-\Psi_m|$. Barron \cite{Barron} formulates
 precise conditions on these
almost minimal values in order to obtain equation (\ref{approxg}),
but we do not
use these conditions here.
The important point we use is that  optimising only one
sigmoidal function each time, Jones' sequence is able to achieve
$O(1/m)$ approximation.
This permits to avoid a global optimisation of all the coefficients
involved nonlinearly.
The linearly involved coefficients $\alpha_m$ and $B_m$
can be computed by
the least square method and the nonlinearly involved
coefficients $A_m$ (which are
two-dimensional like $q$) and $\eta_m$ can be  determined by a
random method.
Such approach allows us to
approximate functions (\ref{fMeinh}) and (\ref{gMeinh}).
The numerical
parameters were $\alpha=8$, $\alpha_1=1$, $\beta_1=0$, $\kappa_1=2$,
$\beta_2=1$, $\kappa_2=0$. In this case system (\ref{RD-1})--(\ref{gMeinh})
with Neumann boundary conditions has an homogeneous equilibrium solution
$u_0=\beta_2/\kappa_1$, $v_0=(\kappa_1^2+\alpha_1\beta_2^2)/(\alpha\beta_2)$.
On the segment $[0,L]$ with $L=60$, and with
$d_u=1$, $d_v=50$, this equilibrium is unstable with respect to some
non-homogeneous perturbations. Using an initial perturbation on $u$ at
the left side ($u(t=0,x)=2u_0$ for $x\in[0,L/10]$), one obtains a
non-homogeneous stationary solution (so-called Turing structure).
Moreover, the perturbation spreads to the right like a wave. See
\cite{Meinhardt2} p. 30. This behaviour is presented in fig. \ref{Meintarget}.
We approximated Meinhardt's model by a gene network, with $600$ genes
($300$ to approximate function (\ref{fMeinh}) and $300$ for (\ref{gMeinh})).
The behaviour is qualitatively similar to the solution of Meinhardt's model.
See fig. \ref{Meinapp}.

The second point we illustrate is the universal pattern
generation problem: given a spatio-temporal pattern, to find a gene
network generating pattern.

We first generate spatio-temporal patterns with one spatial dimension.
The corresponding gene network can be defined as follows:
\begin{equation}\label{sys1D1}
\frac{\partial y_1}{\partial t}=\kappa(\sigma(\theta_1(x))-y_1),
\end{equation}
\begin{equation}\label{sys1D2}
\frac{\partial \bar y_1}{\partial t}=2\kappa(\sigma(\theta_1(x))-\bar y_1),
\end{equation}

\begin{equation}\label{sys1D3}
\frac{\partial u_j}{\partial t}=\lambda(R_j
\sigma(K_jy_1+\bar K_j\bar y_1-\eta_j)-u_j),\ j=1,\dots,m,
\end{equation}
\begin{equation}\label{sys1D4}
\frac{\partial y_{\mathrm{out}}}{\partial t}=
\lambda(\sigma(\sum_{j=1}^mu_j)-y_{\mathrm{out}}).
\end{equation}
Notice that diffusion is absent and
positional information is provided by the spa\-ce-de\-pendent threshold
$\theta_1(x)$.
In this one-dimensional case, the condition on $\theta_1$ means
that this function is strictly monotone.

As it is shown in Sect. ~\ref{sec-prog},
 $t$ and $x$ are functions
of $(y_1,\bar y_1)$. Namely,
\begin{equation}
t=-\frac{1}{\kappa}\log(\frac{\bar y_1}{y_1}-1)
\end{equation}
and
\begin{equation}
x=
\theta_1^{-1}(\sigma^{-1}(\frac{y_1^2}{2y_1-\bar y_1})).
\end{equation}
Thus any $z(t,x)$ can be presented as a function of $(y_1,\bar y_1)$:
\begin{equation}
z(t,x)=z(\mathrm{time}(y_1,\bar y_1),\mathrm{space}(y_1,\bar y_1))=
Z(y_1,\bar y_1).
\end{equation}
Notice that  $y_{\mathrm{out}}$ approximates 
$r=\sigma(\sum_{j=1}^m u_j)$ as $\lambda \to \infty$.
In turn,  $r$ approximates
$\sigma(\sum_{j=1}^mR_j\ 
\sigma(T_jy_1+\bar T_j\bar y_1-\theta_j))$.
Hence, to solve the pattern generation
problem,
we have to determine the coefficients $R_j$, $T_j$, $\bar T_j$, $\theta_j$
such that $\sum_{j=1}^mR_j\sigma(T_jy_1+\bar T_j\bar y_1-\theta_j)$
approximates
$\sigma^{-1}(Z(y_1,\bar y_1))$. It is possible if $Z$
is a continuous function.

The problem is thus to approximate a function of two variables ($y_1$ and
$\bar y_1$).
This problem is intractable with the least square and the random
methods, and we use here
again Jones' iterative approximation method
(see above).


To avoid
singularities at
the lines $y_1=\bar y_1$ and $2y_1=\bar y_1$,
we have approximated this function $Z$ in the image of
the bounded rectangle
  $[T_0,T_1]\times[x_0,x_1]$ by the map
$\displaystyle (t,x)\mapsto (\sigma(\theta_1(x))(1-e^{-\kappa t}),
\sigma(\theta_1(x))(1-e^{-2\kappa t}))$, which is a one-to-one map
 of $(t,x)$.

\begin{figure}[!hb]    
\resizebox{\textwidth}{!}{
\includegraphics{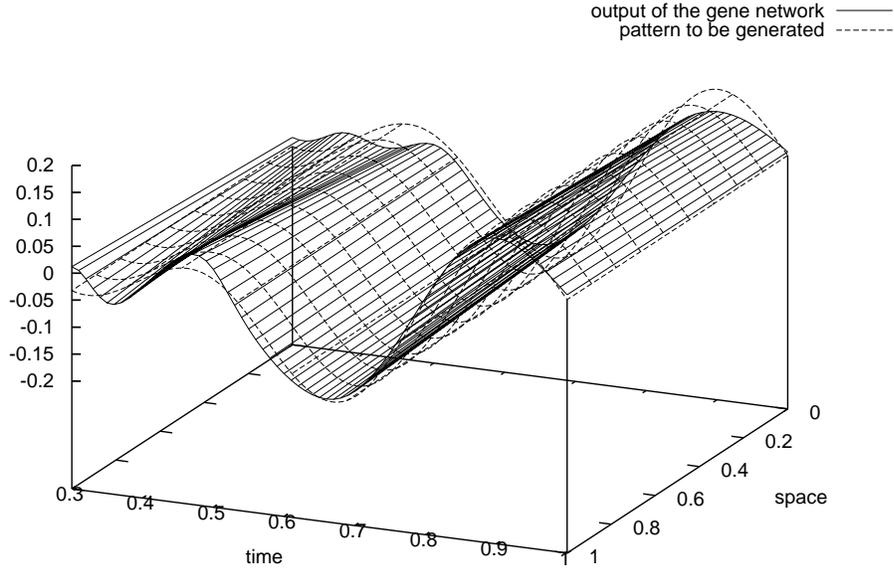}
}
\caption{Generation of $0.1(\sin(8t)+\sin(16t))$ by a gene network.}
\label{unD1}
\end{figure}

\begin{figure}[!hb]
\resizebox{\textwidth}{!}{
\includegraphics{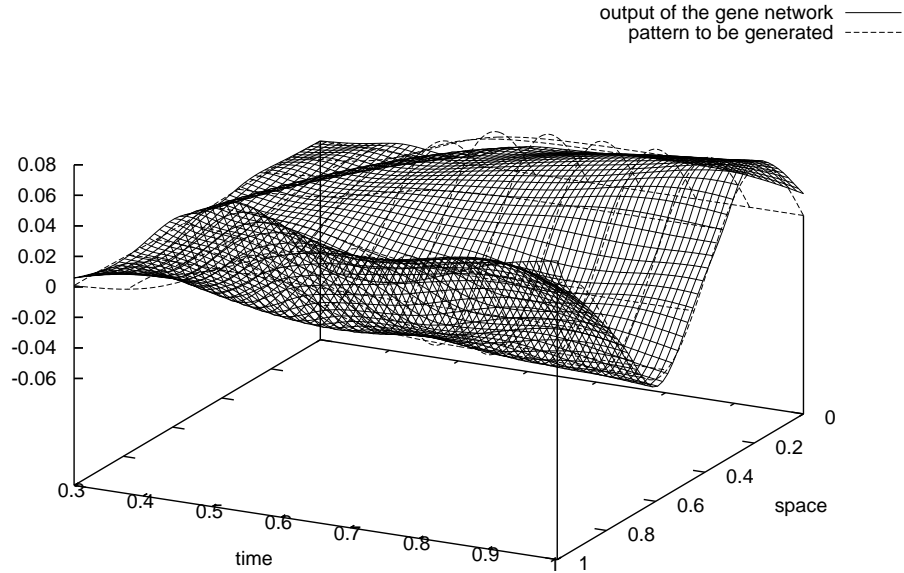}
}
\caption{Generation of $0.025(1+\tanh(10t-0.5))\sin(8x)$ by a gene network.}
\label{unD2}
\end{figure}

Fig. \ref{unD1} and \ref{unD2}  present the output of system
(\ref{sys1D1})--(\ref{sys1D3}) approximating the function
$0.1(\sin(8t)+\sin(16t))$ and $0.025(1+\tanh(10t-0.5))\sin(8x)$,
respectively, for $t\in[0,1]$ and
$x\in[0,1]$. We have used $1000$ sigmoidal functions for these
simulations.

Also we have generated spatio-temporal patterns with $2$ space dimensions.
The corresponding gene circuit is
\begin{equation}\label{sys2D1}
\frac{\partial y_1}{\partial t}=\kappa(\sigma(\theta_1(x))-y_1),
\end{equation}

\begin{equation}\label{sys2D2}
\frac{\partial y_2}{\partial t}=\kappa(\sigma(\theta_2(x))-y_2),
\end{equation}
\begin{equation}\label{sys2D3}
\frac{\partial \bar y_1}{\partial t}=2\kappa(\sigma(\theta_1(x))-\bar y_1),
\end{equation}

\begin{equation}\label{sys2D4}
\frac{\partial u_j}{\partial t}=\lambda(R_j
\sigma(K_j^1y_1+K_j^2y_2+\bar K_j\bar y_1-\eta_j)-u_j),\
j=1,\dots,m,
\end{equation}
\begin{equation}\label{sys2D5}
\frac{\partial y_{\mathrm{out}}}{\partial t}=
\lambda(\sigma(\sum_{j=1}^mu_j)-y_{\mathrm{out}}).
\end{equation}

The time $t$, the  spatial coordinates $x_1$ and
$x_2$ can be expressed as functions of $(y_1,y_2,\bar y_1)$:
\begin{equation}
t=-\frac{1}{\kappa}
\log(\frac{\bar y_1}{y_1}-1),
\end{equation}

\begin{equation}
x_1=
g_1(\sigma^{-1}(\frac{y_1^2}{2y_1-\bar y_1}),
\sigma^{-1}(\frac{y_1y_2}{2y_1-\bar y_1}))
\end{equation}
and
\begin{equation}
x_2=
g_2(\sigma^{-1}(\frac{y_1^2}{2y_1-\bar y_1}),
\sigma^{-1}(\frac{y_1y_2}{2y_1-\bar y_1})).
\end{equation}

Hence, any continuous function $z(t,x_1,x_2)$ can be represented
as a function
of $(y_1,y_2,\bar y_1)$, which has to be approximated by Jones' method in
order
to solve the pattern generation problem. Since $t$,
$x_1$ and
$x_2$ are singular in $y_1=\bar y_1$ and $2y_1=\bar y_1$,
these functions were approximated in the
image of the cubic domain
$[T_0,T_1]\times[x_{1,0},x_{1,1}]\times[x_{2,0},x_{2,1}]$ by  the  map
$(t,x_1,x_2)\mapsto
(\sigma(\theta_1(x))(1-e^{-\kappa t}),\sigma(\theta_2(x))(1-e^{-\kappa t}),
\sigma(\theta_1(x))(1-e^{-2\kappa t}))$.

\begin{figure}[!hb]
\resizebox{\textwidth}{!}{
\includegraphics{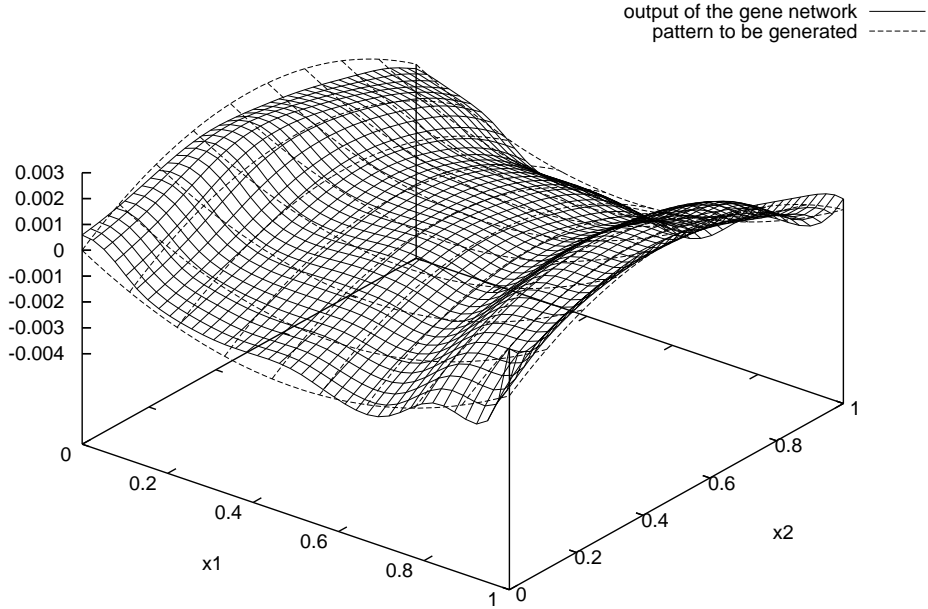}
}
\caption{Generation of $0.01((x_1-0.5)^2-(x_2-0.5)^2)$ by a gene network.}
\label{deuxD}
\end{figure}

Fig. \ref{deuxD}  presents the output of system
(\ref{sys2D1})--(\ref{sys2D5}) approximating the function
$0.01((x_1-0.5)^2-(x_2-0.5)^2)$
for $(x_1,x_2)\in[0,1]^2$. This function is independent of time, but
time-dependent functions have also been approximated (it is not shown).
We have used $1000$ sigmoidal functions for this
simulation.

The last point we illustrate is the superposition principle and its
relation with the conditional complexity (see Sect. ~\ref{sec-prog}).
The superposition Theorem ~\ref{superp} states that a given network
generating a pattern $u(t,x)$ and a given continuous function $F$,
one can device a new network generating $F(u)(t,x)$.
The number of the genes involved in this new network
depends on the complexity of the target pattern.
This complexity can be defined
by the Fourier transform of the pattern \cite{Barron}.
We define the conditional complexity $\mathrm{Comp}(F(u)(t,x)|u(t,x))$ as the
complexity of $F(u)(t,x)$ {\it considered as a function of $u(t,x)$}.
The point is that $\mathrm{Comp}(F(u)(t,x)|u(t,x))$ can be much less than
$\mathrm{Comp}(F(u)(t,x))$. So generating $F(u)(t,x)$ through $u(t,x)$
we may use much less genes than generating $F(u)(t,x)$ directly
(or, if the same gene number is involved, a better precision may be
achieved).

We illustrate this fact by generating $\cos(8t)$ for $t\in[0,2\pi]$.
We produce this time function directly and, moreover, we first
generate $\cos(t)$, then $2\cos^2(t)-1=\cos(2t)$, later
$2\cos^2(2t)-1=\cos(4t)$
and finally $2\cos^2(4t)-1=\cos(8t)$.
The network generating $\cos(8t)$ directly is

\begin{equation}\label{dir1}
\frac{\partial y_1}{\partial t}=1,
\end{equation}

\begin{equation}\label{dir2}
\frac{\partial u_j}{\partial t}=\lambda(R_j\sigma(K_jy_1-\eta_i)-u_j),\
j=1,\dots,m,
\end{equation}

\begin{equation}\label{dir3}
\frac{\partial y_{\mathrm{out}}}{\partial t}=\lambda(\sum_{j=1}^mu_j-
y_{\mathrm{out}}),
\end{equation}
where $R_j$, $K_j$ and $\eta_j$ are chosen so that
$\sum_{j=1}^mR_j\sigma(K_jt-\eta_j)$ approximates $\cos(8t)$.

The network producing $\cos(8t)$ indirectly is
\begin{equation}\label{indir1}
\frac{\partial y_1}{\partial t}=1,
\end{equation}
\begin{equation}
\frac{\partial u_j^1}{\partial t}=
\lambda(R_j^1\sigma(K_j^1y_1-\eta_j^1)-u_j^1),\
j=1,\dots,m_1,
\end{equation}
\begin{equation}
\frac{\partial y_2}{\partial t}=\lambda(\sum_{j=1}^{m_1}u_j^1-y_2),
\end{equation}

\begin{equation}
\frac{\partial u_j^2}{\partial t}=
\lambda(R_j^2\sigma(K_j^2y_2-\eta_j^2)-u_j^2),\
j=1,\dots,m_2,
\end{equation}
\begin{equation}
\frac{\partial y_3}{\partial t}=\lambda(\sum_{j=1}^{m_2}u_j^2-y_3),
\end{equation}

\begin{equation}
\frac{\partial u_j^3}{\partial t}=
\lambda(R_j^2\sigma(K_j^2y_3-\eta_j^2)-u_j^3),\
j=1,\dots,m_2,
\end{equation}
\begin{equation}
\frac{\partial y_4}{\partial t}=\lambda(\sum_{j=1}^{m_2}u_j^3-y_4),
\end{equation}
\begin{equation}
\frac{\partial u_j^4}{\partial t}=
\lambda(R_j^2\sigma(K_j^2y_4-\eta_j^2)-u_j^4),\
j=1,\dots,m_2,
\end{equation}
\begin{equation}\label{indir9}
\frac{\partial y_{\mathrm{out}}}{\partial t}=
\lambda(\sum_{j=1}^{m_2}u_j^4-y_{\mathrm{out}}),
\end{equation}
where $R_j^1$, $K_j^1$ and $\eta_j^1$ have been chosen so that
$\sum_{j=1}^{m_1}R_j^1\sigma(K_j^1t-\eta_j^1)$ approximates $\cos(t)$
and $R_j^2$, $K_j^2$ and $\eta_j^2$ have been chosen so that
$\sum_{j=1}^{m_2}R_j^2\sigma(K_j^2x-\eta_j^2)$  approximates $2x^2-1$.

We have used the same number of equations in the two cases, namely $52$,
(so, $m=50$ in equation (\ref{dir2}))
and we have compared
the precision achieved. For the indirect approximation, we have chose
$m_1=32$ and hence, $m_2=5$.

The target pattern is defined by a function of a single
variable. In this case and with a small number of the genes,
using of Jones' method
is not obligatory.
Actually, here by the least square method for the linear coefficients
and a random choice for the nonlinear ones we achieve a  better precision.

\begin{figure}[!hb]
\resizebox{\textwidth}{!}{
\includegraphics{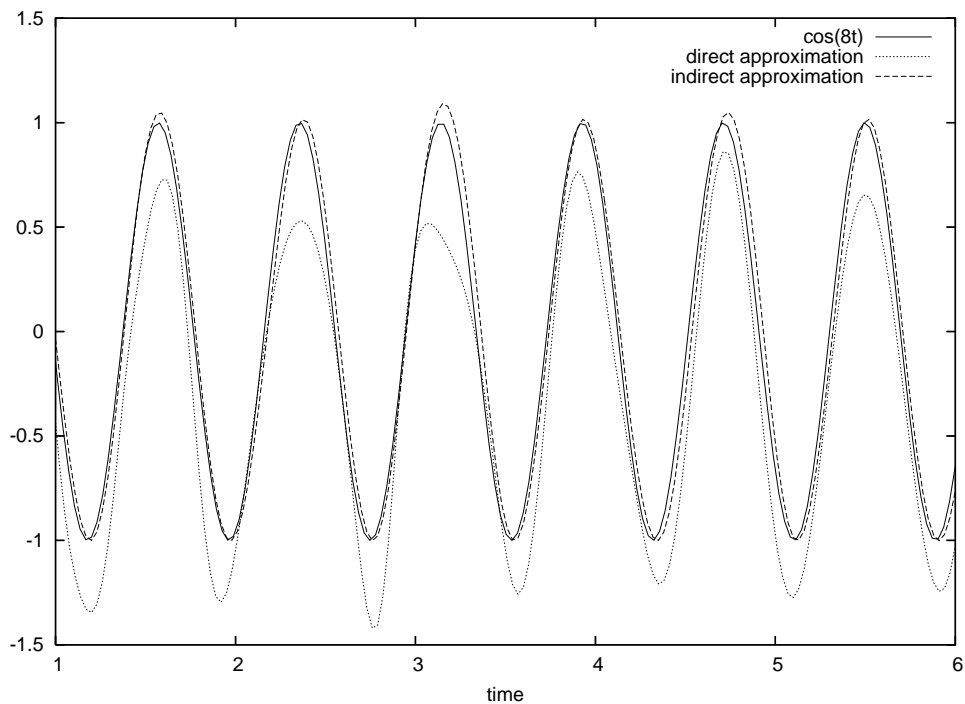}
}
\caption{Improvement of approximation by superposition principle.}
\label{accel}
\end{figure}

Fig. \ref{accel}  presents the results. The patterns computed by systems
(\ref{dir1})--(\ref{dir3}) and (\ref{indir1})--(\ref{indir9})
are denoted respectively "direct approximation" and
"indirect approximation".

\section{Conclusion}

\label{sec-ccl}

It is shown that the genetic networks with binary interaction
of the genes   have  a formidable patterning capacity.
They can produce any spatio-tem\-po\-ral patterns.
Moreover, it is proved that  any reaction-diffusion systems
can be approximated by  genetic circuits. This result
 allows to connect earlier phenomenological mathematical reaction-diffusion
models
 and more biologically realistic genetic circuits.

Let us emphasise that, by these circuits,   pattern programming
can be performed.
This means that, for a given pattern,
 a circuit that builds this pattern,
can be found by  effective and universal algorithms.
One of the most
astonishing biological revelations of the past twenty years is that much
of the basic machinery of development is essentially the same, not
in all vertebrates but in all the major phyla of invertebrates
too \cite{Wolpert}.  We show therefore that this machinery
can be described by simple gene circuit models.

This
pattern programming  holds on
a basic biological principle: on modular organisation of genes.
Genes are organised in blocs.
 Notice that the modular principle is confirmed by experimental data
of molecular biology (see \cite{Alberts,Hop2,Wolpert} and references
therein).

We have demonstrated that this modular structure entails
an interesting  property, which can be named
"superposition principle".
This superposition property means
 that new patterns can always be obtained
by previous (old) patterns.

 As an elementary  example explaining a biological
interpretation of  superposition
principle we can consider  flappers, wings and legs
of tetrapodes. It is well known that  they consist of the same basic
elements (numerus, cubitus, radius,
carpe) but jointed in different ways (see \cite{Rid}).
 Different joinings give wings for birds,  legs for  dogs,
flappers for whales etc.
When mammals penetrated
in water, evolution did not invented flappers from zero.
Evolution used earlier created patterns to obtain
flappers.
Thus, the superposition principle allows us to understand why
the gene number grows relatively slow in evolution
(remind that {\sl Drosophila} has the $14000$ genes,
{\sl C. elegans} has the $19000$ and {\sl Homo sapiens}
has the $30000$ genes).
Indeed, as we have explained in  Sect.  ~\ref{sec-prog} and 
~\ref{sec-num},
the gene number growth  is not directly proportional
to the pattern complexity; this number  is proportional to
conditional pattern complexity relatively already
stored patterns. This conditional quantity  may be small even
if the target pattern is very complex.

We can thus conclude that the modular organisation  of gene interaction
leads to a minimisation of time and genes in a process of invention of new
biological structures.
A famous basic
evolutionary law asserts
that the ontogenesis summarises the philogenesis \cite{Rid}.
 The mathematical results
of this paper suggest that
this law is a direct consequence of gene network organisation.


\appendix

\section{Appendix}

\subsection{Derivation of equations (\ref{u}) and (\ref{v})}

Equations (\ref{gen-RD2-1}) and (\ref{gen-RD2-2}) can be rewritten as
\begin{equation}
\frac{\partial y_i}{\partial t}=\sigma(a_i u  + \gamma_i v  -\theta_i)
  +  d_1 \Delta y_i,
\label{y}
\end{equation}

\begin{equation}
\frac{\partial z_i}
{\partial t}=\sigma(\bar a_i v  + \bar \gamma_i u     -\bar \theta_i)
  +  d_2 \Delta z_i.
\label{z}
\end{equation}

Multiplying the $i$-th equation in
(\ref{y}) by $b_i$ and taking the sum over $i$,
we obtain (\ref{u}) and (\ref{v}).

\subsection{Boundedness of solutions of the Meinhardt
equations}

Solutions of
(\ref{RD-1})-(\ref{RD-2}) stay bounded: they lie in the domain
$D_{C_1, C_2}$ for all times if the corresponding initial data
are in  this domain.

Let us show that condition (\ref{bound}) holds
with appropriate $C_1, C_2$.
Let us choose such $C_i$ that

\begin{equation}
 \alpha C_2 (C_1^2(1+ \alpha_1 C_1^2)^{-1} +\beta_1) < \kappa_1 C_1
\label{C1}
\end{equation}

and
\begin{equation}
  \beta_2 <  (\alpha \beta_1 + \kappa_2)C_2.
\label{C2}
\end{equation}

First we choose a  large
$C_2$ to satisfy (\ref{C2}) and then  we can take a
 constant $C_1$
 large enough to satisfy (\ref{C1}).
 Now we can prove
that the domain $D_{C_1, C_2}$ is an invariant
rectangle for (\ref{RD-1})-(\ref{RD-2}).  This means that on the boundaries
$u=C_1$
and $v=C_2$ the vector field $(f,g)$ is directed inside  $D_{C_1, C_2}$.

This assertion follows from
 (\ref{C1}) and
(\ref{C2}).

\subsection{Proof of Superposition Theorem}

We consider a large network involving a number of genes.
First, it involves the genes $y^s_i, \ i=1,\dots, m_s$,
 where $s=1,\dots,p$, participating in  given networks.
The corresponding dynamics is defined by the equations
\begin{equation}
 \frac{\partial y_i^s(x,t)}{\partial t}=
R_i^s\sigma( \sum_{j=1}^{m_s}
 K_{ij}^s y_j^s - \theta_{i}(x) - \eta_i^s) -\lambda_i^s
y_i^s,\quad x \in \Omega, \quad t \ge 0,
\label{space-threshold1}
\end{equation}
where $s=1,\dots,p$.

Moreover, the large network includes
additional genes $v_k$.
The time evolution of the corresponding
concentrations
$v_k(x,t)$  is defined by the following equations
\begin{equation}
 \frac{\partial v_k(x,t)}{\partial t}=
b_k\sigma(\sum_{j=1}^{p}
M_{kj} y_1^j(x,t) -  \eta_k) - \lambda v_k,\quad x \in \Omega, t > 0.
\label{space-threshold2}
\end{equation}

At last, the genes $v_k(x,t)$ determine the time evolution of the
output gene $y_1$  as follows:
\begin{equation}
 \frac{\partial y_1(x,t)}{\partial t}=
R \sigma(\sum_{j=1}^{m_0}
S_{k} v_k(x,t)) -
 \lambda y_1,\quad x \in \Omega, \quad t > 0.
\label{space-threshold3}
\end{equation}

We set the zero initial conditions for all
the concentrations
$$
    v_k(x,0)=y_j^k(x,0)=y_1(x,0)=0.
$$

Let us prove the auxiliary lemma.

\begin{lemma}
Given a function $z(t) \in C^1[0, T]$,
positive numbers $\epsilon$ and $\delta < T$, there are a function
$w \in C[0, T]$ and a positive coefficient $\lambda$ such that
the solution of the Cauchy problem
\begin{equation}
 \frac{dX(t)}{dt}=
 - \lambda X(t) + w(t), \quad X(0)=0, \quad t \in [0, T]
\label{space-threshold4}
\end{equation}
satisfies
 the following inequality

\begin{equation}
 | X(t) -z(t)| < \epsilon,
  \quad t \in [\delta, T].
\label{space-threshold5}
\end{equation}
\end{lemma}

\debproof
The proof of this lemma is elementary. Indeed,
let us set
$$
    w=\frac{dz(t)}{dt}+
  \lambda z(t), \quad X=z + \tilde X.
$$

Then (\ref{space-threshold4}) entails
$$
\frac{d \tilde X(t)}{dt}=
 - \lambda\tilde X(t), \quad  \tilde X(0)=-z(0).
$$

Thus,
$$
|\tilde X(t)| \le |z(0)|\exp(-\lambda t).
$$

To complete the proof of the lemma, we
 set
$$
\lambda > -\delta^{-1} \log(|z(0)|^{-1}\epsilon).
$$
\finproof

Notice that the lemma also holds
for $z \in C[0, T]$, since any continuous function can
be approximated by a smooth function. Moreover,
if given $z$ is a superposition of the form
$z=z(y^1(t), \dots, y^p(t))$, where $y^s$
are defined by some system of  autonomous differential
equations, then $w$ can also be represented as a superposition:
$w=w(y^1(t), \dots, y^p(t))$.

To finish the proof of   Theorem ~\ref{superp}, it is
sufficient to prove that for any continuous
function of the form $w(x,t)=
w(y_1^1(x,t),  \dots, y_1^p(x,t))$,
where $x \in \Omega,\ t \in [0, T]$ and
$\epsilon > 0$, there exists  such a choice of
the parameters $M_{kj}, \eta_k,
 \lambda, b_k, S_k$ and $R$ in
(\ref{space-threshold2}) and (\ref{space-threshold3})
that the solutions $v_k(x,t)$
of (\ref{space-threshold2}) satisfy
the estimate
\begin{equation}
 | w(y_1^1(x,t),  \dots, y_1^p(x,t)) - R\sigma(\sum_{k=1}^{m_0} S_k v_k)|
< \epsilon,
\label{space-threshold6}
\end{equation}

for any  $x \in \Omega, \ t \in [\delta, T]$.
Using the monotonicity of
$\sigma$ and choosing a sufficiently large $R$, we simplify the last estimate and obtain
\begin{equation}
 | W(y_1^1(x,t),  \dots, y_1^p(x,t)) - \sum_{k=1}^{m_0} S_k v_k| < \epsilon,
  \quad x \in \Omega, \quad t \in [\delta, T],
\label{space-threshold16}
\end{equation}
where $W(x,t)$ is given.

Let  us take sufficiently large $ \lambda > 0$.
Using (\ref{space-threshold16}) and (\ref{space-threshold2}),
we obtain
\begin{equation}
|v_k(x,t) -
 \lambda^{-1}  b_k\sigma(\sum_{j=1}^{p}
M_{kj} y_1^j(x,t) -  \eta_k)| < \epsilon/4.
\label{space-threshold7}
\end{equation}

Denote $\beta_k=\lambda^{-1} S_k b_k$.
Now, to finish the proof, it is sufficient to find  the parameters

$M_{kj}, \eta_k,
 \beta_k$ such that

\begin{equation}
 | W(y_1^1(x,t),  \dots, y_1^p(x,t)) - \sum_{k=1}^{m_0}
  \beta_k\sigma(\sum_{j=1}^p
M_{kj} y_1^j(x,t) -  \eta_k)| < \epsilon/4,
\label{space-threshold8}
\end{equation}
where $x \in \Omega, \ t \in [\delta, T]$.
The existence of this approximation follows
from the multilayered network theory (see Lemma ~\ref{approxlemma}).

\noindent
{\bf Acknowledgments}
\medskip

\noindent
We are thankful to Dr. J. Reinitz (New York) for fruitful
discussion and Dr. V. Volpert (Lyon) for the help.
The paper was supported by program PICS (CNRS and Russia
Academy of Sciences). The second author was supported by the grant
PAST (France).

%

\end{document}